\documentclass[12pt,a4paper,reqno]{amsart}
\usepackage{amssymb,amsmath,amstext,amsthm,graphics}
\usepackage{enumerate}
\usepackage{latexsym}
\usepackage{epsfig}
\usepackage{float}
\usepackage{float}
\usepackage{graphics}
\usepackage{graphicx}
\usepackage{wrapfig}
\usepackage{url}

\theoremstyle{plain} 
\newtheorem{theorem}{Theorem}[section]

\newtheorem{corollary}{Corollary}[section]
\newtheorem{proposition}{Proposition}[section]

\theoremstyle{definition} 

\newcommand{\be}{\begin{eqnarray}}
\newcommand{\ee}{\end{eqnarray}}
\newcommand{\nn}{\nonumber}
\newcommand{\bn}{\begin{eqnarray*}}       
\newcommand{\en}{\end{eqnarray*}}

\newcommand{\Z}{\ensuremath{\mathbb{Z}}}
\newcommand{\R}{\ensuremath{\mathbb{R}}}
\newcommand{\N}{\ensuremath{\mathbb{N}}}

\makeatletter
 \def\@secnumfont{\bfseries}
 \def\section{\@startsection{section}{1}%
   \z@{.7\linespacing\@plus\linespacing}{.5\linespacing}%
   {\normalfont\bfseries\centering}}
\makeatother

\begin{document}

\numberwithin{equation}{section} 

\title[A computer-assisted uniqueness proof for a semilinear elliptic boundary value problem]
	  {A computer-assisted uniqueness proof for a semilinear elliptic boundary value problem}

\author
[P.~J.~McKenna]{Patrick~J.~McKenna}
\author
[F.~Pacella]{Filomena~Pacella}
\author
[M.~Plum]{Michael~Plum}
\author
[D.~Roth]{Dagmar~Roth}

\address{Dep. of Mathematics, 
University of Connecticut, Storrs, 
CT 06269-3009, MSB 328, USA}
\email{mckenna@math.uconn.edu}

\address{Dipartimento di Matematica, 
Universit\`a di Roma ``La Sapienza'', 
P.le A. Moro 2, 00185 Roma, Italy}
\email{pacella@mat.uniroma1.it}

\address{Institut f\"ur Analysis, 
Karlsruher Institut f\"ur Technologie (KIT), 
76128 Karlsruhe, Germany}
\email{michael.plum@kit.edu}

\address{Institut f\"ur Analysis, 
Karlsruher Institut f\"ur Technologie (KIT), 
76128 Karlsruhe, Germany}
\email{dagmar.roth@kit.edu}

\subjclass[2010]{35J25, 35J60, 65N15}

\keywords{semilinear elliptic boundary value problem, uniqueness, computer-assisted proof}

\begin{abstract}
A wide variety of articles, starting with the famous paper \cite{GidNiNir79}, is devoted to the
uniqueness question for the semilinear elliptic boundary value problem $-\Delta u=\lambda u+u^p$ in
$\Omega$, $u>0$ in $\Omega$, $u=0$ on $\partial\Omega$, where $\lambda$ ranges between $0$ and the
first Dirichlet Laplacian eigenvalue. So far, this question was settled in the case of $\Omega$
being a ball and, for more general domains, in the
case $\lambda = 0$. In \cite{McKPacPlR09}, we proposed a
computer-assisted approach to this uniqueness question, which indeed provided a proof in the case
$\Omega=(0,1)^2$, and $p=2$. Due to the high numerical complexity, we were not able in
\cite{McKPacPlR09} to treat higher values of $p$. Here, by a significant reduction of the
complexity, we will prove uniqueness for the case $p=3$.
\end{abstract}

\maketitle

\begin{center} {\itshape Dedicated to the memory of Wolfgang Walter}\end{center}

 \section{Introduction}
 \label{sec:intro}
The semilinear elliptic boundary value problem
\be\label{eq1.1}
- \Delta u = f(u) {\rm ~in~} \Omega, \quad u = 0 {\rm ~on~} \partial \Omega
\ee
has attracted a lot of attention since the 19th century. Questions of existence and multiplicity
have been (are still being) extensively studied by means of variational methods, fixed-point
methods, sub- and supersolutions, index and degree theory, and more.

In this article, we will address the question of {\it uniqueness} of solutions for the more
special problem
\be\label{eq1.2}
\left \{\begin{array}{rcll}
-\Delta u & = & \lambda u+u^p& \text{in } \Omega\\
         u& >& 0& \text{in }\Omega\\
		  u & = & 0 & \text{on }\partial\Omega
        \end{array} \right .
\ee
where $\lambda$ ranges between $0$ and $\lambda_1 (\Omega)$, the first eigenvalue of the Dirichlet
Laplacian. It has been shown in a series of papers \cite{NiNuss85}, \cite{Zhang92}, \cite{Sri93},
\cite{AdiYad94}, \cite{AftPac03} that the solution of \eqref{eq1.2} is indeed unique when
$\Omega$ is a ball, or when $\Omega$ is more general but $\lambda = 0$ (\cite{Zou94},
\cite{Gro00}, \cite{DamGroPac99}, \cite{Dan88}).

We will concentrate on the case where $\Omega = (0,1)^2$ and $p=3$, and prove that uniqueness holds
for the full range $[0,\lambda_1 (\Omega))$ of $\lambda$. Thus, our paper constitutes the first
uniqueness result for this situation. More precisely we prove

\begin{theorem}\label{theo1.1}
 Let $\Omega$ be the unit square in $\R^2$, $\Omega=(0,1)^2$. Then the problem 
\be\label{eq1.3}
\left \{\begin{array}{rcll}-\Delta u & = & \lambda u+u^3& \text{in } \Omega\\
         u& >& 0& \text{in }\Omega\\
		  u & = & 0 & \text{on }\partial\Omega
        \end{array} \right .
\ee
admits only one solution for any $\lambda\in[0,\lambda_1(\Omega))$.
\end{theorem}

{\bf Remark 1.1.} a) A simple scaling argument shows that our uniqueness result carries over to all squares
$\Omega_l:=(0,l)^2$ (and thus, to all squares in $\R^2$): If $u$ is a positive solution of $-\Delta u=\tilde\lambda u+u^3$
in $\Omega_l$, $u=0$ on $\partial\Omega_l$, for some $\tilde\lambda\in[0,\lambda_1(\Omega_l))$, then
$v(x,y):=lu(lx,ly)$ is a solution of \eqref{eq1.3} for $\lambda=\tilde\lambda l^2\in[0,\lambda_1(\Omega))$.
\smallskip

b) Since we also show that the unique solution in the square is nondegenerate, by a result of \cite{Dan88} we deduce that the solution 
is unique also in domains ``close to'' a square.
\smallskip

c) Finally we observe that having shown in \cite{McKPacPlR09} (case  $p=2$) and in this paper (case $p=3$) that the unique solution is 
nondegenerate then uniqueness follows also for other nonlinearities of the type $\lambda u+u^p$ for $p$ close to $2$ and $3$. 
Indeed, by standard arguments (see for example \cite{DamGroPac99}) nonuniqueness of positive solutions in correspondence 
to sequences of exponents converging to $3$ (resp. to  $2$) would imply degeneracy of the solution for $p=3$ (resp. $p=2$). 
\smallskip

Our proof heavily relies on {\it computer-assistance}. Such computer-assisted proofs are receiving
an increasing attention in the recent years since such methods provided results which apparently
could not be obtained by purely analytical means (see
\cite{BreuMcKennaPlum03}, \cite{BreuHorMcKennaPlum06}, \cite{PlumWie02}, \cite{NagNaYa99},
\cite{NaYa95}).

We compute a branch of approximate solutions
and prove existence of a true solution branch close to it, using fixed point techniques. By
eigenvalue enclosure methods, and an additional analytical argument for $\lambda$ close to
$\lambda_1(\Omega)$ we deduce the non-degeneracy of all solutions along this branch, whence
uniqueness follows from the known bifurcation structure of the problem.

In \cite{McKPacPlR09} we give a general description of these computer-assisted means and use them to
obtain the
desired uniqueness result for the case $\Omega = (0,1)^2, p=2$. To make the present paper dealing
with the case $p=3$ more self-contained, we recall parts of the content of \cite{McKPacPlR09} here.
We remark
that the numerical tools used in \cite{McKPacPlR09} turned out not to be sufficient to treat the
case $p=3$.
Now, by some new trick to reduce the numerical complexity, we are able to handle this case.\\

\section{Preliminaries}
\label{sec2}

In the following, let $\Omega = (0,1)^2$. We remark that the results of this section can be carried
over to the more general case of a ``doubly symmetric'' domain; see \cite{McKPacPlR09} for details.
 
First, note that problem \eqref{eq1.2} can equivalently be reformulated as finding a {\it
non-trivial} solution of
\be\label{eq2.1}
\left\{ \begin{array}{rll}
- \Delta u & \hspace{-0,2cm} = \lambda u+|u|^p & {\rm ~in~} \Omega \\
 u & \hspace{-0,2cm} = 0 & {\rm ~on~} \partial \Omega, \end{array} \right. 
\ee
since, for $\lambda<\lambda_1(\Omega)$, by the strong maximum principle (for $- \Delta - \lambda$)
every non-trivial solution of
\eqref{eq2.1} is positive in $\Omega$. In fact, this formulation is better suited for our
computer-assisted approach than \eqref{eq1.2}.

As a consequence of the classical bifurcation theorem of \cite{Rab71} and of the results of
\cite{DamGroPac99} the following result was obtained in
\cite{PacSri03}:

\begin{theorem}\label{theo2.1}
All solutions $u_{\lambda}$ of \eqref{eq1.2} lie on a simple continuous curve $\Gamma$ in
$[0,\lambda_1 (\Omega)) \times C^{1,\alpha} (\bar\Omega)$ joining $(\lambda_1 (\Omega), 0)$ with
$(0, u_0)$, where $u_0$ is the unique solution of \eqref{eq1.2} for $\lambda =0$.
\end{theorem}

We recall that the uniqueness of the solution of \eqref{eq1.2} for $\lambda = 0$ was proved in
\cite{Dan88} and \cite{DamGroPac99}. As a consequence of the previous theorem we have

\begin{corollary}\label{coro2.1}
If all solutions on the curve $\Gamma$ are nondegenerate then problem \eqref{eq1.2} admits only one
solution for every $\lambda \in [0,\lambda_1(\Omega))$.
\end{corollary}

{\bf Proof.} The nondegeneracy of the solutions implies, by the Implicit Function Theorem, that
neither turning points nor secondary bifurcations can exist along $\Gamma$. Then, for every $\lambda
\in [0,\lambda_1 (\Omega))$ there exists only one solution of \eqref{eq1.2} on $\Gamma$. By
Theorem \ref{theo2.1} all solutions are on $\Gamma$, hence uniqueness follows. \hfill $\qed$\\

Theorem \ref{theo2.1} and Corollary \ref{coro2.1} indicate that to prove the uniqueness of the
solution of problem \eqref{eq1.2} for every $\lambda \in [0,\lambda_1 (\Omega))$ it is enough to
construct a branch of nondegenerate solutions which connects $(0,u_0)$ to $(\lambda_1 (\Omega),0)$.
This is what we will do numerically in the next sections with a rigorous computer-assisted proof.\\

However, establishing the nondegeneracy of solutions $u_{\lambda}$ for $\lambda$ close to $\lambda_1
(\Omega)$ numerically can be difficult, due to the fact that the only solution at $\lambda =
\lambda_1 (\Omega)$, which is the identically zero solution, is obviously degenerate because its
linearized operator is $L_0 = - \Delta - \lambda_1$ which has the first eigenvalue equal to zero.
The next proposition shows that there exists a
computable number $\bar{\lambda} (\Omega) \in (0,\lambda_1 (\Omega))$ such that for any $\lambda \in
[\bar{\lambda} (\Omega), \lambda_1 (\Omega))$ problem \eqref{eq1.2} has only one solution which
is also nondegenerate. Of course, from the well-known results of Crandall and Rabinowitz,
\cite{CR1,CR2}, one can establish that for $\lambda$ ``close to'' $\lambda_1$, all solutions $u_{
\lambda }$ are nondegenerate. However, in order to complete our program, we need to calculate a
precise and {\it explicit} estimate of how close they need to be. This allows us to carry out the
numerical computation only in the interval $[0,\bar{\lambda} (\Omega)]$ as we will do later.\\

Let us denote by $\lambda_1 = \lambda_1 (\Omega)$ and $\lambda_2 = \lambda_2 (\Omega)$ the first
and second eigenvalue of the operator $-\Delta$ in $\Omega$ with homogeneous Dirichlet boundary
conditions. We have 
\begin{proposition}\label{pro2.1}
If there exists $\bar{\lambda} \in (0,\lambda_1)$ and a solution $u_{\bar{\lambda}}$ of
\eqref{eq1.2} with $\lambda = \bar{\lambda}$ such that
\be\label{eq2.2}
\| u_{\bar{\lambda}} \|_{\infty} < \left( \frac{\lambda_2 - \lambda_1}{p} \right)^{\frac{1}{p-1}}
\cdot \left( \frac{\bar{\lambda}}{\lambda_1} \right)^{\frac{1}{p-1}}
\ee
then
\be\label{eq2.3}
\| u_{\lambda} \|_{\infty} < \left( \frac{\lambda_2 - \lambda_1}{p} \right)^{\frac{1}{p-1}},
\ee
 and $u_\lambda$ is non-degenerate, for all solutions $u_{\lambda}$ of \eqref{eq1.2} belonging to
the branch $\Gamma_2 \subset \Gamma$
which connects $(\bar{\lambda}, u_{\bar{\lambda}}) $ to $(\lambda_1, 0)$. 
\end{proposition}
(Recall that $\Gamma$ is the unique continuous branch of solutions given by Theorem
\ref{theo2.1}.)\medskip

{\bf Proof.} see \cite{McKPacPlR09}

\begin{corollary}\label{coro2.2}
If on the branch $\Gamma$ there exists a solution $u_{\bar{\lambda}},\ \bar{\lambda} \in (0,\lambda_1)$ such that:
\begin{enumerate}
\item[i)] on the sub-branch $\Gamma_1$ connecting $(0,u_0)$ with $(\bar{\lambda}, u_{\bar{\lambda}})$ all solutions are nondegenerate\\[2ex]
\hspace*{-1,3cm}and
\be\label{eq2.4}
\hspace*{-3,9cm} {\it ii)} \hspace{3cm} \| u_{\bar{\lambda}} \|_{\infty} < \left( \frac{\lambda_2 - \lambda_1}{p} \right)^{\frac{1}{p-1}} \cdot \left( \frac{\bar{\lambda}}{\lambda_1} \right)^{\frac{1}{p-1}},
\ee
\end{enumerate}
then all solutions of \eqref{eq1.2} are nondegenerate, for all $\lambda \in (0,\lambda_1)$, and
therefore problem \eqref{eq1.2} admits only one solution for every $\lambda \in [0,\lambda_1
(\Omega))$.
\end{corollary}

{\bf Proof.} We set $\Gamma = \Gamma_1\cup \Gamma_2$  with $\Gamma_1$ connecting $(0,u_0)$ to
$(\bar\lambda,u_{\bar\lambda})$. On
$\Gamma_1$ we have that all solutions are nondegenerate by i). On the other hand the hypothesis ii)
allows us to apply Proposition \ref{pro2.1} 
which shows nondegeneracy of all solutions on $\Gamma_2$.
Hence there is nondegeneracy all along $\Gamma$ so
the assertion follows from Corollary \ref{coro2.1}. \hfill $\qed$\\

The last corollary suggests the method of proving the uniqueness through computer assistance: first
we construct a branch of nondegenerate ``true'' solutions near approximate ones in a certain
interval $[0,\bar{\lambda}]$ and then verify ii) for the solution $u_{\bar{\lambda}}$. Note that the
estimate \eqref{eq2.4} depends only on $p$ and on the eigenvalues $\lambda_1$ and $\lambda_2$ of the
operator $-\Delta$ in the domain $\Omega$. So the constant on the right-hand side is easily
computable. When $\Omega$ is the unit square which is the case analyzed in the next sections, the
estimate \eqref{eq2.4} becomes:
\bn
\| u_{\bar{\lambda}} \|_{\infty} < \left( \frac{3 \pi^2}{p} \right)^{\frac{1}{p-1}} \cdot \left(
\frac{\bar{\lambda}}{2 \pi^2} \right)^{\frac{1}{p-1}}=\left (\frac{3\bar\lambda}{2p}\right
)^{\frac1{p-1}}
\en
because $\lambda_1 = 2 \pi^2$ and $\lambda_2 = 5 \pi^2$.\\

Fixing $p=3$ we finally get the condition
\be\label{eq2.5}
\| u_{\bar{\lambda}} \|_{\infty} < \sqrt{ \frac{\bar{\lambda}}{2} }.
\ee

\section{The basic existence and enclosure theorem}
\label{sec3}

We start the computer-assisted part of our proof with a basic theorem on existence, local
uniqueness, and non-degeneracy of solutions to problem \eqref{eq2.1}, assuming $p = 3$ now for
simplicity of presentation. In this section, the parameter $\lambda \in [0,\lambda_1 (\Omega))$ is
fixed.

Let $H_0^1 (\Omega)$ be endowed with the inner product $\langle u,v \rangle_{H_0^1} := \langle
\nabla u, \nabla v \rangle_{L^2} + \sigma \langle u, v \rangle_{L^2}$; actually we choose $\sigma = 1$
in this paper, but different (usually positive) choices of $\sigma$ are advantageous or even mandatory
in other applications, whence we keep $\sigma$ as a parameter in the following. Let $H^{-1} (\Omega)$ denote the (topological)
dual of $H_0^1 (\Omega)$, endowed with the usual operator $\sup {\rm -norm}$.

Suppose that an {\it approximate} solution $\omega_{\lambda} \in H_0^1 (\Omega)$ of
problem \eqref{eq2.1} has been computed by numerical means, and that a bound $\delta_{\lambda} > 0$
for its {\it defect} is known, i.e.
\be\label{eq3.1}
\| - \Delta \omega_{\lambda} - \lambda \omega_{\lambda} - | \omega_{\lambda} |^3 \|_{H^{-1}} \le
\delta_{\lambda},
\ee
as well as a constant $K_{\lambda}$ such that
\be\label{eq3.2}
\| v \|_{H_0^1} \le K_{\lambda} \| L_{(\lambda, \omega_{\lambda})} [v] \|_{H^{-1}} {\rm ~~for~all~}
v \in H_0^1 (\Omega).
\ee
Here, $L_{(\lambda, \omega_{\lambda})}$ denotes the operator {\it linearizing} problem
\eqref{eq2.1} at $\omega_{\lambda}$; more generally, for $(\lambda, u) \in \R \times H_0^1
(\Omega)$, let the linear operator $L_{(\lambda,u)} : H_0^1 (\Omega) \to H^{-1} (\Omega)$ be defined
by 
\be\label{eq3.3}
L_{(\lambda,u)} [v] := - \Delta v - \lambda v - 3 |u| uv \qquad(v \in H_0^1 (\Omega)).
\ee
The practical computation of bounds $\delta_{\lambda}$ and $K_{\lambda}$ will be addressed in
Sections 6, 7 and 8.

Let $C_{4}$ denote a norm bound (embedding constant) for the embedding $H_0^1 (\Omega)
\hookrightarrow L^{4} (\Omega)$, which is bounded since $\Omega \subset \R^2$. $C_{4}$ can be
calculated e.g. according to the explicit formula given in \cite[Lemma 2]{Plum08}. Finally, let
\bn
\gamma:= 3 C^3_{4}\; .
\en
In our example case where $\Omega = (0,1)^2$, the above-mentioned explicit formula gives (with the
choice $\sigma:=1$)
\bn
\gamma = \frac{3\sqrt2}{4\left( \pi^2 + 1\right)^{3/4} } \left( < \frac{1}{5}
\right) \; .
\en

\begin{theorem}\label{theo3.1} Suppose that some $\alpha_{\lambda} > 0$ exists such that
\be\label{eq3.4}
\delta_{\lambda} \le \frac{\alpha_{\lambda}}{K_{\lambda}}- \gamma \alpha_{\lambda}^2 \left( \left\|
\omega_{\lambda} \right\|_{L^{4}} + C_{4} \alpha_{\lambda} \right) 
\ee
and
\be\label{eq3.5}
2 K_{\lambda} \gamma \alpha_{\lambda} \left( \| \omega_{\lambda} \|_{L^{4}} + C_{4}
\alpha_{\lambda} \right) < 1 \; .
\ee

Then, the following statements hold true:
\begin{enumerate}
\item[a)] (existence) There exists a solution $u_{\lambda} \in H_0^1 (\Omega)$ of problem
\eqref{eq2.1} such that
\be\label{eq3.6}
\| u_{\lambda} - \omega_{\lambda} \|_{H_0^1} \le \alpha_{\lambda} \; .
\ee
\item[b)] (local uniqueness) Let $\eta > 0$ be chosen such that \eqref{eq3.5} holds with
$\alpha_{\lambda} + \eta$ instead of $\alpha_{\lambda}$. Then,
\be\label{eq3.7}
\left. \begin{array}{l}  u \in H_0^1 (\Omega) {\rm ~solution~of~} \eqref{eq2.1} \\[0.5ex]
                  \| u - \omega_{\lambda}\|_{H_0^1} \le \alpha_{\lambda} + \eta \end{array} \right\} \Longrightarrow u = u_{\lambda} \; .
\ee
\item[c)] (nondegeneracy)
\be\label{eq3.8}
\left. \begin{array}{l}  u \in H_0^1 (\Omega)  \\[0.5ex]
                  \| u - \omega_{\lambda}\|_{H_0^1} \le \alpha_{\lambda} \end{array} \right\}
\Longrightarrow L_{(\lambda, u)} : H_0^1 (\Omega) \to H^{-1} (\Omega) {\rm ~~is~bijective},
\ee
whence in particular $L_{(\lambda, u_{\lambda})}$ is bijective (by \eqref{eq3.6}).
\end{enumerate}
\end{theorem}
For a proof, see \cite{McKPacPlR09}.
\begin{corollary}\label{coro3.1} Suppose that \eqref{eq3.4} and \eqref{eq3.5} hold, and in addition
that $\| \omega_{\lambda} \|_{H_0^1} > \alpha_{\lambda}$. Then, the solution $u_{\lambda}$ given by
Theorem 3.1 is non-trivial (and hence positive).
\end{corollary}

{\bf Remark 3.1.}\label{rem3.1} 
a) The function $\psi (\alpha):= \frac{\alpha}{K_{\lambda}} - \gamma \alpha^2 (\| \omega_{\lambda}
\|_{L^{4}} + C_{4} \alpha)$ has obviously a positive maximum at $\bar\alpha = \frac 1{3C_4}\left
(\sqrt{\|\omega_\lambda\|_{L^4}^2+\frac{3C_4}{K_\lambda\gamma}}-\|\omega_\lambda\|_{L^4}
\right )$, and the crucial condition \eqref{eq3.4} requires that
\begin{multline}
\label{eq3.9}
\delta_{\lambda} \le \psi (\bar\alpha) = \frac {4C_4+\gamma K_\lambda\|\omega_\lambda\|_{L^4}^2}{
K_\lambda\left (\sqrt{\gamma K_\lambda(\gamma K_\lambda\|\omega_\lambda\|_{L^4}^2+3 C_4)}+\gamma
K_\lambda\|\omega_\lambda\|_{L^4}\right )}\cdot \\ \frac 1{
\left (\sqrt{\gamma K_\lambda(\gamma K_\lambda\|\omega_\lambda\|_{L^4}^2+3 C_4)}+\gamma
K_\lambda\|\omega_\lambda\|_{L^4}+6C_4\right ) }\ \ ,\\
\end{multline}

i.e. $\delta_{\lambda}$ has to be sufficiently small. According to \eqref{eq3.1}, this means that
$\omega_{\lambda}$ must be computed with sufficient accuracy, which leaves the ``hard work'' to the
computer!

Furthermore, a ``small'' defect bound $\delta_{\lambda}$ allows (via \eqref{eq3.4}) a ``small''
error bound $\alpha_{\lambda}$, if $K_{\lambda}$ is not too large.\smallskip

b) If moreover we choose the {\it minimal} $\alpha_{\lambda}$ satisfying \eqref{eq3.4}, then the
additional condition \eqref{eq3.5} follows automatically, which can be seen as follows:
the minimal choice of $\alpha_{\lambda}$ shows that $\alpha_{\lambda} \leq
\bar{\alpha}$. We have
\begin{multline}
2 K_{\lambda} \gamma \bar{\alpha} (\| \omega_{\lambda} \|_{L^{4}} + C_{4}
\bar{\alpha})=\\1-\frac {C_4}{3C_4+2\gamma K_\lambda \|\omega_\lambda\|_{L^4}^2 +2\sqrt{\gamma
K_\lambda(\gamma K_\lambda\|\omega_\lambda\|_{L^4}^2+3C_4)}\|\omega_\lambda\|_{L^4}}<1
\end{multline}
and thus condition \eqref{eq3.5} is satisfied.

Since we will anyway try to find $\alpha_{\lambda}$ (satisfying \eqref{eq3.4}) close to the minimal
one, condition \eqref{eq3.5} is ``practically'' always satisfied if \eqref{eq3.4} holds.
(Nevertheless, it must of course be checked.) \smallskip

\section{The branch $(u_{\lambda})$}
\label{sec4}

Fixing some $\bar{\lambda} \in (0,\lambda_1 (\Omega))$ (the actual choice of which is made on the
basis of Proposition \ref{pro2.1}; see also Section 5), we assume now
that for {\it every} $\lambda \in [0,\bar{\lambda}]$ an approximate solution $\omega_{\lambda} \in
H_0^1 (\Omega)$ is at hand, as well as a defect bound $\delta_{\lambda}$ satisfying \eqref{eq3.1},
and a bound $K_{\lambda}$ satisfying \eqref{eq3.2}. Furthermore, we assume now that, for every
$\lambda \in [0,\bar{\lambda}]$, some $\alpha_{\lambda} > 0$ satisfies \eqref{eq3.4} and
\eqref{eq3.5}, and the additional non-triviality condition $\| \omega_{\lambda} \|_{H_0^1} >
\alpha_{\lambda}$ (see Corollary 3.1). We suppose that some {\it uniform} ($\lambda$-independent)
$\eta > 0$ can be chosen such that \eqref{eq3.5} holds with $\alpha_{\lambda} + \eta$ instead of
$\alpha_{\lambda}$ (compare Theorem \ref{theo3.1} b)). Hence Theorem \ref{theo3.1} gives a positive
solution $u_{\lambda} \in H_0^1 (\Omega)$ of problem \eqref{eq2.1} with the properties
\eqref{eq3.6}, \eqref{eq3.7}, \eqref{eq3.8}, for every $\lambda \in [0,\bar{\lambda}] $. 

Finally, we assume that the approximate solution branch $([0,\bar{\lambda}] \to H_0^1 (\Omega)$,
$\lambda \mapsto \omega_{\lambda})$ is continuous, and that $([0,\bar{\lambda}] \to \R, ~ \lambda
\mapsto \alpha_{\lambda})$ is lower semi-continuous.

In Sections 6, 7 and 8, we will address the actual computation of such
branches $(\omega_{\lambda}), (\delta_{\lambda}), (K_{\lambda}), $ $(\alpha_{\lambda})$.

So far we  know nothing about continuity or smoothness of $([0,\bar{\lambda}] \to H_0^1 (\Omega), ~
\lambda \mapsto u_{\lambda})$, which however we will need to conclude that $(u_{\lambda})_{\lambda
\in [0,\bar{\lambda}]}$ coincides with the sub-branch $\Gamma_1$ introduced in Corollary
\ref{coro2.2}.

\begin{theorem}\label{theo4.1}
The solution branch
\bn
\left\{ \begin{array}{ccc}
[0,\bar{\lambda}] &\to &H_0^1 (\Omega) \\ \lambda &\mapsto &u_{\lambda} \end{array} \right\}
\en
is continuously differentiable.
\end{theorem}

{\it Proof:} The mapping
\bn
\mathcal{F} : \left\{ \begin{array}{ccc}
\R  \times H_0^1 (\Omega)&\to &H^{-1} (\Omega) \\( \lambda, u)  &\mapsto &- \Delta u - \lambda u -
|u|^3  \end{array} \right\}
\en
is continuously differentiable, with $(\partial \mathcal{F} / \partial u)(\lambda, u) =
L_{(\lambda, u)}$ (see \eqref{eq3.3}), and $\mathcal{F} (\lambda, u_{\lambda}) = 0$ for all $\lambda
\in [0,\bar{\lambda}]$. Using the Mean Value Theorem one can show that $L_{(\lambda, u)}$
depends indeed continuously on $(\lambda, u)$; see \cite[Lemma 3.1]{McKPacPlR09} for details.

It suffices to prove the asserted smoothness locally. Thus, fix $\lambda_0 \in [0,\bar{\lambda}]$.
Since $L_{(\lambda_0, u_{\lambda_0})}$ is bijective by Theorem \ref{theo3.1} c), the Implicit
Function Theorem gives a $C^1$-smooth solution branch
\bn
\left\{ \begin{array}{ccc}
(\lambda_0 - \varepsilon, \lambda_0 + \varepsilon) &\to &H_0^1 (\Omega) \\ \lambda  &\mapsto & \hat{u}_{\lambda}  \end{array} \right\}
\en
to problem \eqref{eq2.1}, with $\hat{u}_{\lambda_0} = u_{\lambda_0}$. By \eqref{eq3.6}, 
\be\label{eq4.1}
\| \hat{u}_{\lambda_0} - \omega_{\lambda_0} \|_{H_0^1} \le \alpha_{\lambda_0} \; .
\ee

Since $\hat{u}_{\lambda}$ and $\omega_{\lambda}$ depend continuously on $\lambda$, and
$\alpha_{\lambda}$ lower semi-continuously, \eqref{eq4.1} implies
\bn
\| \hat{u}_{\lambda} - \omega_{\lambda} \|_{H_0^1} \le \alpha_{\lambda} + \eta \; \; \; \;(\lambda
\in [0,\bar{\lambda}] \cap (\lambda_0 - \tilde{\varepsilon}, \lambda_0 + \tilde{\varepsilon}))
\en
for some $\tilde{\varepsilon} \in (0,\varepsilon)$. Hence Theorem \ref{theo3.1} b) provides
\bn
\hat{u}_{\lambda} = u_{\lambda} \; \; \; \;(\lambda \in [0,\bar{\lambda]} \cap (\lambda_0 - \tilde{\varepsilon}, \lambda_0 + \tilde{\varepsilon})) \; ,
\en
implying the desired smoothness in some neighborhood of $\lambda_0$ (which of course is one-sided
if $\lambda_0 = 0$ or $\lambda_0 = \bar{\lambda})$. \hfill $\qed$ 
\smallskip

As a consequence of Theorem \ref{theo4.1}, $(u_{\lambda})_{\lambda \in [0,\bar{\lambda}]}$ is a
continuous solution curve connecting the point $(0,u_0)$ with $(\bar{\lambda}, u_{\bar{\lambda}})$,
and thus must coincide with the sub-branch $\Gamma_1$, connecting these two points, of the unique
simple continuous curve $\Gamma$ given by Theorem \ref{theo2.1}. Using Theorem \ref{theo3.1} c), we
obtain

\begin{corollary}\label{coro4.1} On the sub-branch $\Gamma_1$ of $\Gamma$ which connects $(0,u_0)$ with $(\bar{\lambda}, u_{\bar{\lambda}})$, all solutions are nondegenerate.
\end{corollary}

Thus, if we can choose $\bar{\lambda}$ such that condition \eqref{eq2.4} holds true, Corollary
\ref{coro2.2} will give the desired uniqueness result.

\section{Choice of $\bar{\lambda}$}
\label{sec5}

We have to choose $\bar{\lambda}$ such that condition \eqref{eq2.4} is satisfied. For this purpose,
we use computer-assistance again. With $x_M$ denoting the intersection of the symmetry axes of the
(doubly symmetric) domain $\Omega$, i.e. $x_M = \left( \frac{1}{2}, \frac{1}{2} \right)$ for $\Omega
= (0,1)^2$, we choose $\bar{\lambda} \in (0,\lambda_1 (\Omega))$, not too close to $\lambda_1
(\Omega)$, such that our {\it approximate} solution $\omega_{\bar{\lambda}}$ satisfies
\be\label{eq5.1}
\omega_{\bar{\lambda}} (x_M) < \left( \frac{\lambda_2 (\Omega) - \lambda_1 (\Omega) }{3}
\right)^{\frac{1}{2}} \cdot \left( \frac{\bar{\lambda}}{\lambda_1 (\Omega)}
\right)^{\frac{1}{2}} \; ,
\ee
with ``not too small'' difference between right- and left-hand side. Such a $\bar{\lambda}$ can be found within a few numerical trials.

Here, we impose the additional requirement
\be\label{eq5.2}
\omega_{\bar{\lambda}} \in H^2 (\Omega) \cap H_0^1 (\Omega) \; ,
\ee
which is in fact a condition on the numerical method used to compute $\omega_{\bar{\lambda}}$.
(Actually, condition \eqref{eq5.2} could be avoided if we were willing to accept additional
technical effort.) Moreover, exceeding \eqref{eq3.1}, we will now need an $L^2$-bound
$\hat{\delta}_{\bar{\lambda}}$ for the defect:
\be\label{eq5.3}
\| - \Delta \omega_{\bar{\lambda}} - \bar{\lambda} \omega_{\bar{\lambda}} - | \omega_{\bar{\lambda}}
|^3 \|_{L^2} \le \hat{\delta}_{\bar{\lambda}} \; .
\ee
Finally, we note that $\Omega$ is {\it convex}, and hence in particular
$H^2$-regular, whence every solution $u \in H_0^1 (\Omega)$ of problem \eqref{eq2.1} is in $H^2
(\Omega)$. 

Using the method described in Section 3, we obtain, by Theorem \ref{theo3.1} a), a positive
solution $u_{\bar{\lambda}} \in H^2 (\Omega) \cap H_0^1 (\Omega)$ of problem \eqref{eq2.1}
satisfying
\be\label{eq5.4}
\| u_{\bar{\lambda}} - \omega_{\bar{\lambda}} \|_{H_0^1} \le \alpha_{\bar{\lambda}} \; ,
\ee
provided that \eqref{eq3.4} and \eqref{eq3.5} hold, and that $\| \omega_{\bar{\lambda}} \|_{H_0^1} >
\alpha_{\bar{\lambda}}$.

Now we make use of the explicit version of the Sobolev embedding $H^2 (\Omega) \hookrightarrow C
(\bar{\Omega})$ given in \cite{Plum92}. There, explicit constants $\hat{C}_0$, $\hat{C}_1$,
$\hat{C_2}$ are computed such that
\bn
\| u \|_{\infty} \le \hat{C}_0 \| u \|_{L^2} + \hat{C}_1 \| \nabla u \|_{L^2} + \hat{C}_2 \| u_{xx}
\|_{L^2} \; \; \; \; {\rm for~all~} u \in H^2 (\Omega) \; ,
\en
with $\| u_{xx} \|_{L^2}$ denoting the $L^2$-Frobenius norm of the Hessian matrix $u_{xx}$. E.g.
for $\Omega = (0,1)^2$, \cite{Plum92} gives
\bn
\hat{C}_0 = 1, ~ \hat{C}_1 = 1.1548 \cdot \sqrt{\frac{2}{3}} \le 0.9429, ~ \hat{C}_2 = 0.22361
\cdot \sqrt{\frac{28}{45}} \le 0.1764  \; .
\en
Moreover, $\| u_{xx} \|_{L^2} \le \| \Delta u\|_{L^2}$ for $u \in H^2 (\Omega) \cap H_0^1 (\Omega)$
since $\Omega$ is convex (see e.g. \cite{Lad68}). Consequently, 
 
\be\label{eq5.5}
\| u_{\bar{\lambda}} - \omega_{\bar{\lambda}} \|_{\infty} \le \hat{C}_0 \| u_{\bar{\lambda}}
-\omega_{\bar{\lambda}} \|_{L^2} + \hat{C}_1 \| u_{\bar{\lambda}} - \omega_{\bar{\lambda}}
\|_{H_0^1} +  \hat{C}_2 \| \Delta u_{\bar{\lambda}} - \Delta \omega_{\bar{\lambda}} \|_{L^2} .
\ee
To bound the last term on the right-hand side, we first note that
\be\label{eq5.6}
\| | u_{\bar{\lambda}} |^3 - |\omega_{\bar{\lambda}} |^3 \|_{L^2} &=& \bigg\| 3 \int_0^1 |
\omega_{\bar{\lambda}}  + t (u_{\bar{\lambda}}  - \omega_{\bar{\lambda}} ) |
(\omega_{\bar{\lambda}}  + t(u_{\bar{\lambda}} - \omega_{\bar{\lambda}} )) dt \cdot
(u_{\bar{\lambda}} - \omega_{\bar{\lambda}} ) \bigg\|_{L^2} \nn \\[1ex]
&\le& 3 \int_0^1 \| | \omega_{\bar{\lambda}} + t(u_{\bar{\lambda}}  - \omega_{\bar{\lambda}}
)|^{2} \cdot | u_{\bar{\lambda}}  - \omega_{\bar{\lambda}} | \|_{L^2} dt \nn  \\[1ex]
&\le& 3 \int_0^1 \|  \omega_{\bar{\lambda}} + t(u_{\bar{\lambda}}  - \omega_{\bar{\lambda}})
\|^{2}_{L^{6}} \| u_{\bar{\lambda}}  - \omega_{\bar{\lambda}}  \|_{L^{6}} dt \nn \\[1ex]
&\le& 3 \int_0^1 \left( \| \omega_{\bar{\lambda}}  \|_{L^{6}} + t C_{6} \alpha_{\bar{\lambda}}
\right)^{2}  dt  \cdot C_{6} \alpha_{\bar{\lambda}}\\[1ex]
& = &3C_6\left (\|\omega_\lambda\|_{L^6}^2+C_6\|\omega_\lambda\|_{L^6}\alpha_{\bar\lambda}+\frac
13C_6^2 \alpha_{\bar\lambda}^2\right )\alpha_{\bar\lambda} ,
\ee
using \eqref{eq5.4} and an embedding constant $C_{6}$ for the embedding $H_0^1(\Omega)
\hookrightarrow L^{6} (\Omega)$ in the last but one line; see e.g. \cite[Lemma 2]{Plum08} for its
computation. Moreover, by \eqref{eq2.1} and \eqref{eq5.3}, 
\be\label{eq5.7}
\| \Delta u_{\bar{\lambda}} - \Delta \omega_{\bar{\lambda}}   \|_{L^2} \le
\hat{\delta}_{\bar{\lambda}}   + \bar{\lambda} \| u_{\bar{\lambda}}   - \omega_{\bar{\lambda}} 
\|_{L^2} + \| | u_{\bar{\lambda}}  |^3 - |\omega_{\bar{\lambda}}  |^3 \|_{L^2} .
\ee
Using \eqref{eq5.4} - \eqref{eq5.7}, and the Poincar\'e inequality
\be\label{eq5.8}
\|u\|_{L^2} \le \frac{1}{\sqrt{\lambda_1 (\Omega) + \sigma}} \| u \|_{H_0^1} \qquad (u \in H_0^1 (\Omega)), 
\ee
we finally obtain

\begin{multline}\label{eq5.9}
\| u_{\bar{\lambda}}  -\omega_{\bar{\lambda}}  \|_{\infty} \le \\ \left[ \frac{\hat{C}_0 +
\bar{\lambda} \hat{C}_2 }{\sqrt{\lambda_1 (\Omega) + \sigma}}  + \hat{C}_1 + 3 C_{6} \hat{C}_2 
 \left (\|\omega_\lambda\|_{L^6}^2+C_6\|\omega_\lambda\|_{L^6}\alpha_{\bar\lambda}+\frac 13
C_6^2\alpha_{\bar\lambda}^2\right )\right ] \cdot \alpha_{\bar{\lambda}} + \hat{C}_2
\hat{\delta}_{\bar{\lambda}} ,
\end{multline}
and the right-hand side is ``small'' if $\alpha_{\bar{\lambda}} $ and $\hat{\delta}_{\bar{\lambda}} 
$ are ``small'', which can (again) be achieved by sufficiently accurate numerical computations.\\

Finally, since
\[
u_{\bar{\lambda}}   (x_M) \le \omega_{\bar{\lambda}}  (x_M) + \| u_{\bar{\lambda}}  - \omega_{\bar{\lambda}}   \|_{\infty} , 
\]
\eqref{eq5.9} yields an upper bound for $u_{\bar{\lambda}} (x_M)$ which is ``not too much'' larger
than $\omega_{\bar{\lambda}}   (x_M)$. Hence, since
$u_{\bar\lambda}(x_M)=\|u_{\bar\lambda}\|_\infty$ by \cite{GidNiNir79}, condition \eqref{eq2.4} can
easily be checked, and \eqref{eq5.1} (with ``not too small'' difference between right- and left-hand
side) implies a good chance that this check will be successful; otherwise, $\bar{\lambda}$ has to be
chosen a bit larger.

\section{Computation of $\omega_{\lambda}$ and $\delta_{\lambda}$ for fixed $\lambda$}
\label{sec6}
 
In this  section we report on the computation of an approximate solution $\omega_{\lambda}  \in H^2
(\Omega) \cap H_0^1 (\Omega)$ to problem \eqref{eq2.1}, and of bounds $\delta_{\lambda} $ and
$K_{\lambda} $ satisfying \eqref{eq3.1} and \eqref{eq3.2}, where $\lambda \in [0,\lambda_1
(\Omega))$ is {\it fixed} (or one of {\it finitely} many values). We will again restrict ourselves
to the unit square $\Omega = (0,1)^2$.

An approximation $\omega_{\lambda} $ is computed by a {\it Newton iteration} applied to problem
\eqref{eq2.1}, where the linear boundary value problems
\be\label{eq6.1}
L_{(\lambda, \omega_{\lambda}^{(n)})} [v_n] = \Delta \omega_{\lambda}^{(n)} + \lambda
\omega_{\lambda}^{(n)} + |\omega_{\lambda}^{(n)}|^3
\ee
occurring in the single iteration steps are solved approximately by an ansatz
\be\label{eq6.2}
v_n (x_1,x_2) = \sum\limits_{i,j = 1}^N \alpha_{ij}^{(n)} \sin (i \pi x_1) \sin (j \pi x_2) 
\ee
and a Ritz-Galerkin method (with the basis functions in \eqref{eq6.2}) applied to problem
\eqref{eq6.1}. The update $\omega_{\lambda}^{(n+1)} := \omega_{\lambda}^{(n)} + v_n$ concludes the
iteration step.

The Newton iteration is terminated when the coefficients $\alpha_{ij}^{(n)}$ in \eqref{eq6.2} are
``small enough'', i.e. their modulus is below some pre-assigned tolerance.

To {\it start} the Newton iteration, i.e. to find an appropriate $\omega_{\lambda}^{(0)}$ of the
form \eqref{eq6.2}, we first consider some $\lambda$ close to $\lambda_1 (\Omega)$, and choose
$\omega_{\lambda}^{(0)}(x_1,x_2) = \alpha \sin (\pi x_1) \sin (\pi x_2)$; with an appropriate choice
of
$\alpha > 0$ (to be determined in a few numerical trials), the Newton iteration will ``converge'' to
a non-trivial approximation $\omega^{(\lambda)}$. Then, starting at this value, we diminish
$\lambda$ in small steps until we arrive at $\lambda = 0$, while in each of these steps the
approximation $\omega^{(\lambda)}$ computed in the {\it previous} step is taken as a start of the
Newton iteration. In this way, we find approximations $\omega_{\lambda}$ to problem \eqref{eq2.1}
for ``many'' values of $\lambda$. Note that all approximations $\omega_{\lambda}$ obtained in this
way are of the form \eqref{eq6.2}. 

\vspace{0.75cm}

The computation of an $L^2$-defect bound $\hat{\delta}_{\lambda}$ satisfying
\be\label{eq6.3}
\| - \Delta \omega_{\lambda} - \lambda  \omega_{\lambda} - | \omega_{\lambda} |^3 \|_{L^2} \le
\hat{\delta}_{\lambda} 
\ee
amounts to the computation of an integral over $\Omega$. 

Due to \cite{GidNiNir79} every solution of \eqref{eq2.1} is symmetric with respect to reflection
at the axes $x_1=\frac 12$ and $x_2=\frac 12$. Therefore it is useful to look for approximate
solutions of the form
\be\label{eq6.4}
\omega_\lambda(x_1,x_2) = \sum\limits_{\substack{i,j=1\\i,j \text{ odd}}}^N\alpha_{ij}\sin(i\pi
x_1)\sin(j\pi x_2).
\ee

Using sum formulas for $\sin$ and $\cos$ one obtains for all $n\in\N_0, x\in\R$
\[\sin((2n+1)\pi x)=\left (2\sum\limits_{k=1}^n\cos(2k\pi x)+1\right )\sin(\pi x)\]
and thus $\omega_\lambda$ can be written as follows:
\begin{multline}\label{eq6.5}
\omega_\lambda(x_1,x_2) = \alpha_{11}\sin(\pi x_1)\sin(\pi x_2)+\\\sum\limits_{k,l=1}^{\left
\lfloor{\frac
{N-1}2}\right \rfloor}\alpha_{2k+1,2l+1}\left (2\sum\limits_{i=1}^k \cos(2i\pi x_1)+1\right )\left
(2\sum\limits_{j=1}^l\cos(2j\pi x_2)+1\right )\sin(\pi x_1)\sin(\pi x_2).
\end{multline}
Since $\cos(x)$ ranges in $[-1,1]$ and $\sin(\pi x_1)\sin(\pi x_2)$ is positive
for $(x_1,x_2)\in \Omega=(0,1)^2$, $\omega_\lambda$ will be positive if
\be\label{eq6.6}
\alpha_{11}+\sum\limits_{k,l=1}^{\left \lfloor{\frac
{N-1}2}\right \rfloor}\alpha_{2k+1,2l+1}\left ([-2k+1,2k+1]\right )\left
([-2l+1,2l+1]\right )\subset(0,\infty).
\ee
Condition \eqref{eq6.6} can easily be checked using interval arithmetic and is indeed always
satisfied for our approximate solutions, since $\alpha_{11}$ turns out to be ``dominant'' and the
higher coefficients decay quickly. Hence $\omega_\lambda$ is positive and one can omit the
modulus in the computations. Therefore the integral in \eqref{eq6.3} can be
computed in closed form, since only products of trigonometric functions occur in the integrand.
After calculating them, various sums $\sum_{i=1}^N$ remain to be evaluated. In order to obtain a
{\it rigorous} bound $\hat{\delta}_\lambda$, these computations (in contrast to those for obtaining
$\omega_\lambda$ as described above) need to be carried out in {\it interval} arithmetic
\cite{Kla93,Rump02}, to take rounding errors into account.

Note that the complexity in the evaluation of the defect integral in \eqref{eq6.3}, without any
further modifications, is
$O(N^{12})$ due to the term $\omega_\lambda^3$. Using some trick, it is however possible to reduce
the complexity to $O(N^6)$:

Applying the sum formulas $\sin(a)\sin(b)=\frac 12[\cos(a-b)-\cos(a+b)]$ and
$\cos(a)\cos(b)=\frac 12[\cos(a-b)+\cos(a+b)]$ one obtains:
\begin{multline*}
\sin(i_1\pi x)\sin(i_2\pi x)\sin(i_3\pi x)\sin(i_4\pi x)\sin(i_5\pi x)\sin(i_6\pi
x) = \\
-\frac 1{32}\sum\limits_{\substack{\sigma_2,\sigma_3,\sigma_4,\\\sigma_5,\sigma_6\in\{-1,1\}}}
\sigma_2\sigma_3\sigma_4\sigma_5\sigma_6\cos\left (\left (i_1+\sigma_2 i_2+\sigma_3 i_3+\sigma_4
i_4+\sigma_5 i_5+\sigma_6 i_6\right )\pi x\right ).
\end{multline*}
Since $\int_0^1\cos(n\pi x)\,dx = \left \{\begin{array}{cl}1 & \text{ for } n=0\\ 0 & \text{ for }
n\in\Z\backslash\{0\}\end{array}\right \}=:\delta_n$, we get

\begin{eqnarray*}
 \int\limits_\Omega\omega_\lambda(x_1,x_2)^6\,d(x_1,x_2)& = & \frac
1{1024}\sum\limits_{\sigma_2,\ldots,\sigma_6\in\{-1,1\}}\sum\limits_{\rho_2,\ldots,\rho_6\in\{-1,1
\}}\sigma_2\cdot\ldots\cdot\sigma_6\cdot\rho_2\cdot\ldots\cdot\rho_6\cdot \\
&&\sum\limits_{i_1,\ldots,
i_6=1 } ^N \sum\limits_{j_1,\ldots,j_6=1}^N\delta_{i_1+\sigma_2i_2+\ldots+\sigma_6i_6}\delta_{
j_1+\rho_2j_2+\ldots+\rho_6j_6}\alpha_{i_1j_1}\cdot\ldots\cdot\alpha_{i_6j_6}.
\end{eqnarray*}

Setting $\alpha_{ij}:=0$ for $(i,j)\in\Z^2\backslash\{1,\ldots,N\}^2$ the previous sum can be
rewritten as

\begin{multline*}
 \frac 1{1024}\sum\limits_{\substack{\sigma_2,\ldots,\sigma_6,\\\rho_2\ldots,\rho_6\in\{-1,1\}}}
\sigma_2\cdot\ldots\cdot\sigma_6\cdot\rho_2\cdot\ldots\cdot\rho_6\cdot \\
\sum\limits_{k=-2N+1}^{3N}
\sum\limits_{l=-2N+1}^{3N}\left (\sum\limits_{\substack{i_1+\sigma_2 i_2\\+\sigma_3 i_3 =
k}}\ \sum\limits_{\substack{j_1+\rho_2 j_2\\+\rho_3 j_3=l}}
\alpha_{i_1j_1}\alpha_{i_2j_2}\alpha_{i_3j_3}\right )
\left (\sum\limits_{\substack{\sigma_4 i_4+\sigma_5 i_5\\+\sigma_6 i_6 =
-k}}\ \sum\limits_{\substack{\rho_4j_4+\rho_5 j_5\\+\rho_6 j_6=-l}}
\alpha_{i_4j_4}\alpha_{i_5j_5}\alpha_{i_6j_6}\right ).
\end{multline*}

For fixed $\sigma_i,\rho_i,k$ and $l$ each of the two double-sums in parentheses is $O(N^4)$.
Since they are independent, the product is still $O(N^4)$. The sums over $k$ and $l$ then give
$O(N^6)$, whereas the sums over $\sigma_i$ and $\rho_i$ do not change the complexity. \\
Moreover the sum $\sum_{k=-2N+1}^{3N}$ is only
\[\left \{\begin{array}{cl}\sum_{k=3}^{3N}&\text{ if } \sigma_2 = 1,\sigma_3 = 1\\[1.5ex]
           \sum_{k=2-N}^{2N-1}& \text{ if }\sigma_2\cdot\sigma_3 = -1\\[1.5ex]
			\sum_{k=-2N+1}^{N-2} & \text{ if }\sigma_2 = -1,\sigma_3 = -1.
          \end{array}
\right .\]
Similarly, also certain constellations of $\sigma_4,\sigma_5,\sigma_6$ reduce the
$k$-sum, and of course analogous reductions are possible for the $l$-sum. Since $\alpha_{ij}=0$ if
$i$ or $j$ is even, the result does not
change if the sum is only taken over odd values of $i_n,j_n,k$ and $l$.
\smallskip

{\bf Remark 6.1.}\label{rem6.1}  a) Computing trigonometric sums in an efficient way is an object of investigation since a very long time, but up to our knowledge the above
complexity reduction has not been published before.

b) As an alternative to the closed form integration described above, we also tried quadrature for computing the defect integral,
but due to the necessity of computing a safe remainder term bound in this case, we ended up in a very high numerical effort, since a large number of quadrature
 points had to be chosen. So {\it practically} closed-form integration turned out to be more efficient, although its complexity (as $N\to\infty$) is 
higher than the quadrature complexity.

\smallskip

Once an $L^2$-defect bound $\hat{\delta}_{\lambda}$ (satisfying \eqref{eq6.3}) has been computed, an
$H^{-1}$-defect bound $\delta_{\lambda}$ (satisfying \eqref{eq3.1}) is easily obtained via the
embedding
\be\label{eq6.7}
\| u \|_{H^{-1}} \le \frac{1}{\sqrt{\lambda_1 (\Omega) + \sigma}} \| u \|_{L^2} \qquad (u \in L^2 (\Omega))
\ee
which is a result of the corresponding dual embedding \eqref{eq5.8}. Indeed, \eqref{eq6.3} and
\eqref{eq6.7} imply that
\bn
\delta_{\lambda} := \frac{1}{\sqrt{\lambda_1 (\Omega) + \sigma}} \hat{\delta}_{\lambda} 
\en
satisfies \eqref{eq3.1}.

The estimate \eqref{eq6.7} is suboptimal but, under practical aspects, seems to be the most suitable way for obtaining
an $H^{-1}$-bound for the defect. At this point we also wish to remark that, as an alternative to the weak solutions approach used in this paper,
we could also have aimed at a computer-assisted proof for {\it strong} solutions (see [23]), leading to $H^2$-and $C^0$-error bounds;
in this case an $L^2$-bound is needed directly (rather than an $H^{-1}$-bound).

\section{Computation of $K_\lambda$ for fixed $\lambda$}
\label{sec7}

For computing a constant $K_{\lambda}$ satisfying \eqref{eq3.2}, we use the isometric isomorphism
\be\label{eq7.1}
\Phi : \left\{ \begin{array}{cll}
H_0^1 (\Omega) & \to & H^{-1} (\Omega) \\
u&\mapsto& - \Delta u + \sigma u \end{array} \right\},
\ee
and note that $\Phi^{-1} L_{(\lambda, \omega_{\lambda})} : H_0^1 (\Omega) \to H_0^1 (\Omega)$ is $\langle \cdot,\cdot \rangle_{H_0^1}$-symmetric since
\be\label{eq7.2}
\langle \Phi^{-1} L_{(\lambda,\omega_{\lambda})} [u], v \rangle_{H_0^1} = \int_{\Omega} \left[
\nabla u \cdot \nabla v - \lambda uv - 3 |\omega_{\lambda} | \omega_{\lambda} u v \right] dx,
\ee
and hence selfadjoint. Since $\|
L_{(\lambda,\omega_\lambda)}[u]\|_{H^{-1}}=\|\Phi^{-1}L_{(\lambda,\omega_\lambda)}[u]\|_{H_0^1}$,
\eqref{eq3.2} thus holds for any
\be\label{eq7.3}
K_{\lambda} \ge \left[ \min \left\{ | \mu| : \mu {\rm ~is~in~the~spectrum~of~} \Phi^{-1} L_{(\lambda, \omega_{\lambda})} \right\} \right]^{-1}, 
\ee
provided the min is positive.

A particular consequence of \eqref{eq7.2} is that
\be\label{eq7.4}
\langle \left( I - \Phi^{-1} L_{(\lambda, \omega_{\lambda})} \right) [u], u \rangle_{H_0^1} = \int_{\Omega} W_{\lambda} u^2 dx \qquad (u \in H_0^1 (\Omega))
\ee
where
\be\label{eq7.5}
W_{\lambda} (x) := \sigma + \lambda + 3|\omega_{\lambda} (x) | \omega_{\lambda} (x) \qquad (x
\in  \Omega).
\ee
Note that, due to the positivity of our approximate solutions $\omega_\lambda$ established in
Section 6, the modulus can be omitted here, which again facilitates numerical computations.
Choosing a {\itshape positive} parameter $\sigma$ in the $H_0^1$-product (recall that we actually chose
$\sigma:=1$), we obtain $W_{\lambda} > 0$ on $\bar{\Omega}$.
Thus, \eqref{eq7.4} shows that all eigenvalues $\mu$ of $\Phi^{-1} L_{(\lambda, \omega_{\lambda})}$
are less than $1$, and that its essential spectrum consists of the single point $1$. Therefore,
\eqref{eq7.3} requires the computation of {\it eigenvalue bounds} for the eigenvalue(s) $\mu$
neighboring $0$.\\
Using the transformation $\kappa = 1/(1-\mu)$, the eigenvalue problem $\Phi^{-1} L_{(\lambda,
\omega_{\lambda})} [u] = \mu u$ is easily seen to be equivalent to
\bn
- \Delta u + \sigma u = \kappa W_{\lambda} u,
\en
or, in weak formulation,
\be\label{eq7.6}
\langle u,v \rangle_{H_0^1} = \kappa \int_{\Omega} W_{\lambda} u v dx \qquad (v \in H_0^1 (\Omega)),
\ee
and we are interested in bounds to the eigenvalue(s) $\kappa$ neighboring $1$. It is therefore
sufficient to compute two-sided bounds to the first $m$ eigenvalues $\kappa_1 \le \dots \le
\kappa_m$ of problem \eqref{eq7.6}, where $m$ is (at least) such that $\kappa_m > 1$. In all our
practical examples, the computed enclosures $\kappa_i \in [\underline{\kappa}_i, \bar{\kappa}_i ]$
are such that $\bar{\kappa}_1 < 1 < \underline{\kappa}_2$, whence by \eqref{eq7.3} and $\kappa =
1/(1-\mu)$ we can choose
\be\label{eq7.7}
K_{\lambda} := \max \left\{ \frac{\bar{\kappa_1}}{1 - \bar{\kappa_1}} , \frac{\underline{\kappa}_2}{\underline{\kappa}_2 - 1} \right\}. 
\ee

{\bf Remark 7.1.}
By \cite{GidNiNir79} and the fact that $\omega_\lambda$ is symmetric with respect to reflection at
the axes $x_1=\frac 12$ and $x_2 =\frac 12$, all occurring function spaces can be replaced by
their intersection with the class of reflection symmetric functions. This has the advantage that
some eigenvalues $\kappa_i$ drop out, which possibly reduces the constant $K_\lambda$.
\vspace{0,5cm}

The desired {\it eigenvalue bounds} for problem \eqref{eq7.6} can be obtained by computer-assisted
means of their own. For example, {\it upper} bounds to $\kappa_1, \dots , \kappa_m$ (with $m \in \N$
given) are easily and efficiently computed by the {\it Rayleigh-Ritz} method \cite{Rek80}:\\
Let $\tilde{\varphi}_1, \dots, \tilde{\varphi}_m \in H_0^1 (\Omega)$ denote linearly independent
trial functions, for example approximate eigenfunctions obtained by numerical means, and form the
matrices
\bn
A_1 := ( \langle \tilde{\varphi}_i,\tilde{\varphi}_j \rangle_{H_0^1})_{i,j = 1, \dots, m}, \qquad 
A_0 := \left( \int\limits_\Omega W_{\lambda} \tilde{\varphi}_i \tilde{\varphi}_j \,dx \right)_{i,j =
1, \dots, m}.
\en
Then, with $\Lambda_1 \le \dots \le \Lambda_m$ denoting the eigenvalues of the matrix eigenvalue problem
\bn
A_1 x = \Lambda A_0 x
\en
(which can be enclosed by means of verifying numerical linear algebra; see \cite{Beh87}), the Rayleigh-Ritz method gives
\bn
\kappa_i \le \Lambda_i \text{ for } i = 1,\dots, m.
\en
However, also {\it lower} eigenvalue bounds are needed, which constitute a more complicated task
than upper bounds. The most accurate method for this purpose has been proposed by Lehmann
\cite{Leh63}, and improved by Goerisch concerning its range of applicability \cite{BehGoe94}. Its
numerical core is again (as in the Rayleigh-Ritz method) a matrix eigenvalue problem, but the
accompanying analysis is more involved. In particular, in order to compute lower bounds to the first
$m$ eigenvalues, a {\it rough} lower bound to the $(m+1)$-st eigenvalue must be known already. This
a priori information can usually be obtained via a {\it homotopy method} connecting a simple ``base
problem'' with known eigenvalues to the given eigenvalue problem, such that all eigenvalues increase
(index-wise) along the homotopy; see \cite{Plum97} or \cite{BreuHorMcKennaPlum06} for details on
this method, a detailed description of which would be beyond the scope of this article. In fact,
\cite{BreuHorMcKennaPlum06} contains the newest version of the homotopy method, where only very
small ($2 \times 2$ or even $1 \times 1$) matrix eigenvalue problems need to be treated rigorously
in the course of the homotopy.\\
Finding a base problem for problem \eqref{eq7.6}, and a suitable homotopy connecting them, is rather
simple here since $\Omega$ is a bounded rectangle, whence the eigenvalues of $- \Delta$ on $H_0^1
(\Omega)$ are known: We choose a constant upper bound $c_0$ for $|\omega_{\lambda}|
\omega_{\lambda} =\omega_\lambda^2$ on $\Omega$, and the coefficient homotopy
\bn
W_{\lambda}^{(s)} (x) := \sigma + \lambda + 3 [(1-s) c_0 + s \omega_{\lambda} (x)^2 
] \qquad (x \in \Omega, 0 \le s \le 1).
\en
Then, the family of eigenvalue problems
\[
-\Delta u + \sigma u  = \kappa^{(s)} W_{\lambda}^{(s)} u
\]
connects the explicitly solvable constant-coefficient base problem $(s=0)$ to problem \eqref{eq7.6}
$(s=1)$, and the eigenvalues increase in $s$, since the Rayleigh quotient does, by Poincar\'{e}'s
min-max principle.\\

\section{Computation of branches $(\omega_{\lambda} ), (\delta_{\lambda}), (K_{\lambda}), (\alpha_{\lambda})$}
\label{sec8}

In the previous section we described how to compute approximations $\omega_{\lambda} $ for a grid of
finitely many values of $\lambda$ within $[0,\lambda_1 (\Omega))$. After selecting $\bar{\lambda}$
(among these) according to Section 5, we are left with a grid
\bn
0 = \lambda^0 < \lambda^1 < \cdots < \lambda^M = \bar{\lambda}
\en
and approximate solutions $\omega^i = \omega_{\lambda^i} \in H_0^1 (\Omega) \cap L^{\infty}
(\Omega) ~ (i = 0, \dots, M )$. Furthermore, according to the methods described in the previous
sections, we can compute bounds $\delta^i = \delta_{\lambda^i}$ and $K^i = K_{\lambda^i}$ such that
\eqref{eq3.1} and \eqref{eq3.2} hold at $\lambda = \lambda^i$.

Now we define a piecewise linear (and hence continuous) approximate solution branch
$([0,\bar{\lambda}] \to H_0^1 (\Omega), \lambda \mapsto \omega_{\lambda} )$ by
\be\label{eq8.1}
\omega_{\lambda}  := \frac{ \lambda^i - \lambda}{\lambda^i - \lambda^{i-1}} \omega^{i-1} + 
\frac{ \lambda - \lambda^{i-1}}{\lambda^i - \lambda^{i-1}} \omega^{i} \qquad(\lambda^{i-1} <
\lambda < \lambda^i ,  i = 1, \dots, M).
\ee

To compute corresponding defect bounds $\delta_{\lambda} $, we fix $i \in \{ 1, \dots, M\}$ and
$\lambda \in [\lambda^{i-1}, \lambda^i]$, and let $t:= (\lambda - \lambda^{i-1})/(\lambda^i -
\lambda^{i-1}) \in [0,1]$, whence
\be\label{eq8.2}
\lambda = (1-t) \lambda^{i-1} + t \lambda^i, ~~ \omega_{\lambda}  = (1-t) \omega^{i-1} + t \omega^i.
\ee

Using the classical linear interpolation error bound we obtain, for fixed $x \in \Omega$,
\begin{align}\label{eq8.3}
 \big| \omega_{\lambda} (x)  ^3 & -\left [(1-t) \omega^{i-1} (x)^3 + t \omega^i (x)^3 \right ]
\big| \le
\notag \\ 
&\le  \frac{1}{2} \max_{s \in [0,1]} \left| \frac{d^2}{ds^2} \left [ (1-s) \omega^{i-1} (x) + s
\omega^i (x)\right  ]^3 \right| \cdot t (1-t)  \notag \\
&\le  \frac{3}{4}\max_{s \in [0,1]} \left[  (1-s) \omega^{i-1} (x) + s \omega^i (x)
\right] \cdot (\omega^i (x) - \omega^{i-1} (x))^2  \notag \\
&\le \frac{3}{4} \max \left\{ \| \omega^{i-1}\|_{\infty}, \|  \omega^i\|_{\infty}
\right\} \| \omega^i - \omega^{i-1} \|_{\infty}^2 , 
\end{align}

\begin{align}\label{eq8.4}
 \big| \lambda \omega_{\lambda} (x) & - \left[(1-t) \lambda^{i-1} \omega^{i-1} (x)  + t \lambda^i  \omega^i (x) \right] \big|  \nn \\
&\le \frac{1}{2} \max_{s \in [0,1]} \left| \frac{d^2}{ds^2} \left[ ((1-s) \lambda^{i-1} + s \lambda^i) ((1-s) \omega^{i-1} (x) + s \omega^i (x)) \right] \right| \cdot t (1-t) \nn \\ 
&\le \frac{1}{4} (\lambda^i - \lambda^{i-1}) \| \omega^i - \omega^{i-1} \|_{\infty} .
\end{align}

Since $\| u \|_{H^{-1}} \le C_1 \| u \|_{\infty}$ for all $u \in L^{\infty} (\Omega)$, with $C_1$
denoting an embedding constant for the embedding $H_0^1 (\Omega) \hookrightarrow  L^1 (\Omega)$
(e.\,g. $C_1 = \sqrt{|\Omega |} C_2$), \eqref{eq8.3} and \eqref{eq8.4} imply
\begin{multline}\label{eq8.5}
\| \omega_{\lambda}^3 - [(1-t) (\omega^{i-1})^3 + t(\omega^i )^3 ] \|_{H^{-1}} \\
\le \frac{3}{4}C_1 \max \{ \| \omega^{i-1} \|_{\infty} , \| \omega^i \|_{\infty} \}
\| \omega^i - \omega^{i-1} \|_{\infty}^2 = : \rho_i,
\end{multline}

\be\label{eq8.6}
\| \lambda \omega_{\lambda} - [(1-t) \lambda^{i-1} \omega^{i-1} + t \lambda^i \omega^i ] \|_{H^{-1}} \le
 \frac{1}{4} C_1 (\lambda^i - \lambda^{i-1}) \| \omega^i - \omega^{i-1} \|_{\infty} = : \tau_i.
\ee
Now \eqref{eq8.2}, \eqref{eq8.5}, \eqref{eq8.6} give
\begin{align}\label{eq8.7}
&\| - \Delta \omega_{\lambda}  - \lambda \omega_{\lambda}  -  \omega_{\lambda}  ^3 \|_{H^{-1}} 
\nn \\
&\le (1-t) \| - \Delta \omega^{i-1} - \lambda^{i-1} \omega^{i-1} - (\omega^{i-1} )^3 \|_{H^{-1}} + t
\| - \Delta \omega^{i} - \lambda^{i} \omega^{i} - (\omega^{i} )^3 \|_{H^{-1}} + \tau_i + \rho_i \nn
\\
&\le \max \{ \delta^{i-1} , \delta^i´ \} + \tau_i + \rho_i =: \delta_{\lambda}.
\end{align}
Thus, we obtain a branch $(\delta_{\lambda} )_{\lambda \in [0,\bar{\lambda}]}$ of defect bounds
which is constant on each subinterval $[\lambda^{i-1},\lambda^i]$. In the points $\lambda^1, \dots,
\lambda^{M-1}, \delta_{\lambda} $ is possibly doubly defined by \eqref{eq8.7}, in which case we
choose the smaller of the two values. Hence, $([0,\bar{\lambda}] \to \R, \lambda \mapsto
\delta_{\lambda} )$ is lower semi-continuous.\\
Note that $\delta_{\lambda} $ given by \eqref{eq8.7} is ``small'' if $\delta^{i-1}$ and $\delta^i$
are small (i.e. if the approximations $\omega^{i-1}$ and $\omega^i$ have been computed with
sufficient accuracy; see Remark 3.1a)) {\it and} if $\rho_i, \tau_i$ are small (i.e. if the grid is
chosen sufficiently fine; see \eqref{eq8.5}, \eqref{eq8.6}).

\vspace{1cm}

In order to compute bounds $K_{\lambda} $ satisfying \eqref{eq3.2} for $\lambda \in
[0,\bar{\lambda}]$, with $\omega_{\lambda} $ given by \eqref{eq8.1}, we fix $i \in \{ 1, \dots,
M-1\}$ and $\lambda \in \left[ \frac{1}{2} (\lambda^{i-1} + \lambda^i), \frac{1}{2} (\lambda^i +
\lambda^{i+1}) \right]$. Then,

\be\label{eq8.8}
\begin{split}
| \lambda - \lambda^i | & \le \frac{1}{2} \max \{ \lambda^i - \lambda^{i-1}, \lambda^{i+1} - \lambda^i \} =: \mu_i, \\
\| \omega_{\lambda} - \omega^i \|_{H_0^1} & \le \frac{1}{2} \max \{ \| \omega^i - \omega^{i-1} \|_{H_0^1}, \| \omega^{i+1} - \omega^{i} \|_{H_0^1} \} =: \nu_i, 
\end{split}
\ee

whence a coefficient perturbation result given in \cite[Lemma 3.2]{McKPacPlR09} implies: If

\be\label{eq8.9}
\zeta_i := K^i \left[ \frac{1}{\lambda_1 (\Omega) + \sigma} \mu_i + 2 \gamma (\|
\omega^i\|_{L^{4}} + C_{4} \nu_i )\nu_i \right] < 1,
\ee

then \eqref{eq3.2} holds for
\be\label{eq8.10}
K_{\lambda} := \frac{K^i}{1 - \zeta_i} .
\ee

Note that \eqref{eq8.9} is indeed satisfied if the grid is chosen sufficiently fine, since then
$\mu_i$ and $\nu_i$ are ``small'' by \eqref{eq8.8}.\\
Analogous estimates give $K_{\lambda}$ also on the two remaining half-intervals $[0,\frac{1}{2}
\lambda^1]$ and \linebreak $[\frac{1}{2} (\lambda^{M-1} + \lambda^M), \lambda^M]$.\\
Choosing again the smaller of the two values at the points $\frac{1}{2} (\lambda^{i-1} +
\lambda^i)\ (i = 1, \dots, M)$ where $K_{\lambda}$ is possibly doubly defined by \eqref{eq8.10}, we
obtain a lower semi-continuous, piecewise constant branch $([0,\bar{\lambda}] \rightarrow \R,
\lambda \mapsto K_{\lambda})$.

\vspace{1cm}

According to the above construction, both $\lambda \mapsto \delta_{\lambda}$ and $\lambda \mapsto
K_{\lambda}$ are constant on the $2M$ half-intervals. Moreover, \eqref{eq8.1} implies that, for $i =
1, \dots, M$,
\[
\| \omega_{\lambda} \|_{L^{4}}\le \left\{\begin{array}{l}
\max \{ \| \omega^{i-1} \|_{L^{4}}, \frac{1}{2} (\| \omega^{i-1} \|_{L^{4}} + \| \omega^i
\|_{L^{4}} ) \} {\rm ~for~} \lambda \in [\lambda^{i-1}, \frac{1}{2} (\lambda^{i-1} + \lambda^i)]
\\
\max \{ \frac{1}{2} ( \| \omega^{i-1} \|_{L^{4}}+ \| \omega^{i} \|_{L^{4}}),  \| \omega^i
\|_{L^{4}} ) \} {\rm ~for~} \lambda \in [ \frac{1}{2} ( \lambda^{i-1} + \lambda^i), \lambda^i] 
\end{array}
\right\}
\]
and again we choose the smaller of the two values at the points of double definition.

Using these bounds, the crucial inequalities \eqref{eq3.4} and \eqref{eq3.5} (which have to be
satisfied for all $\lambda \in [0,\bar{\lambda}]$) result in {\it finitely} many inequalities which
can be fulfilled with ``small'' and piecewise constant $\alpha_{\lambda}$ if $\delta_{\lambda}$ is
sufficiently small, i.e. if $\omega^0,\dots, \omega^M$ have been computed with sufficient accuracy
(see Remark 3.1a)) and if the grid has been chosen sufficiently fine (see \eqref{eq8.5} -
\eqref{eq8.7}). Moreover, since $\lambda \mapsto \delta_{\lambda}$, $\lambda \mapsto K_{\lambda}$
and the above piecewise constant upper bound for $\| \omega_{\lambda} \|_{L^{4}}$ are lower
semi-continuous, the structure of the inequalities \eqref{eq3.4} and \eqref{eq3.5} clearly shows
that also $\lambda \mapsto \alpha_{\lambda}$ can be chosen to be lower semi-continuous, as required
in Section 4. Finally, since \eqref{eq3.5} now consists in fact of {\it finitely} many strict
inequalities, a uniform ($\lambda$-independent) $\eta > 0$ can be chosen in Theorem 3.1b), as needed
for Theorem 4.1.

\section{Numerical results}
\label{sec9}

All computations have been performed on an AMD Athlon Dual Core
4800+ (2.4GHz) processor, using MATLAB (version R2010a) and the interval toolbox INTLAB
\cite{Rump02}. For some of the time consuming nested sums occurring in the computations, we used
moreover mexfunctions to outsource these calculations to C++. For these parts of the program we
used C-XSC \cite{Kla93} to verify the results. Our source code can be found on our
webpage\footnote{\url{http://www.math.kit.edu/iana2/~roth/page/publ/en}}.

In the following, we report on some more detailed numerical results.

\begin{figure}[H]
 \includegraphics[angle=90,width=0.8\textwidth,clip]{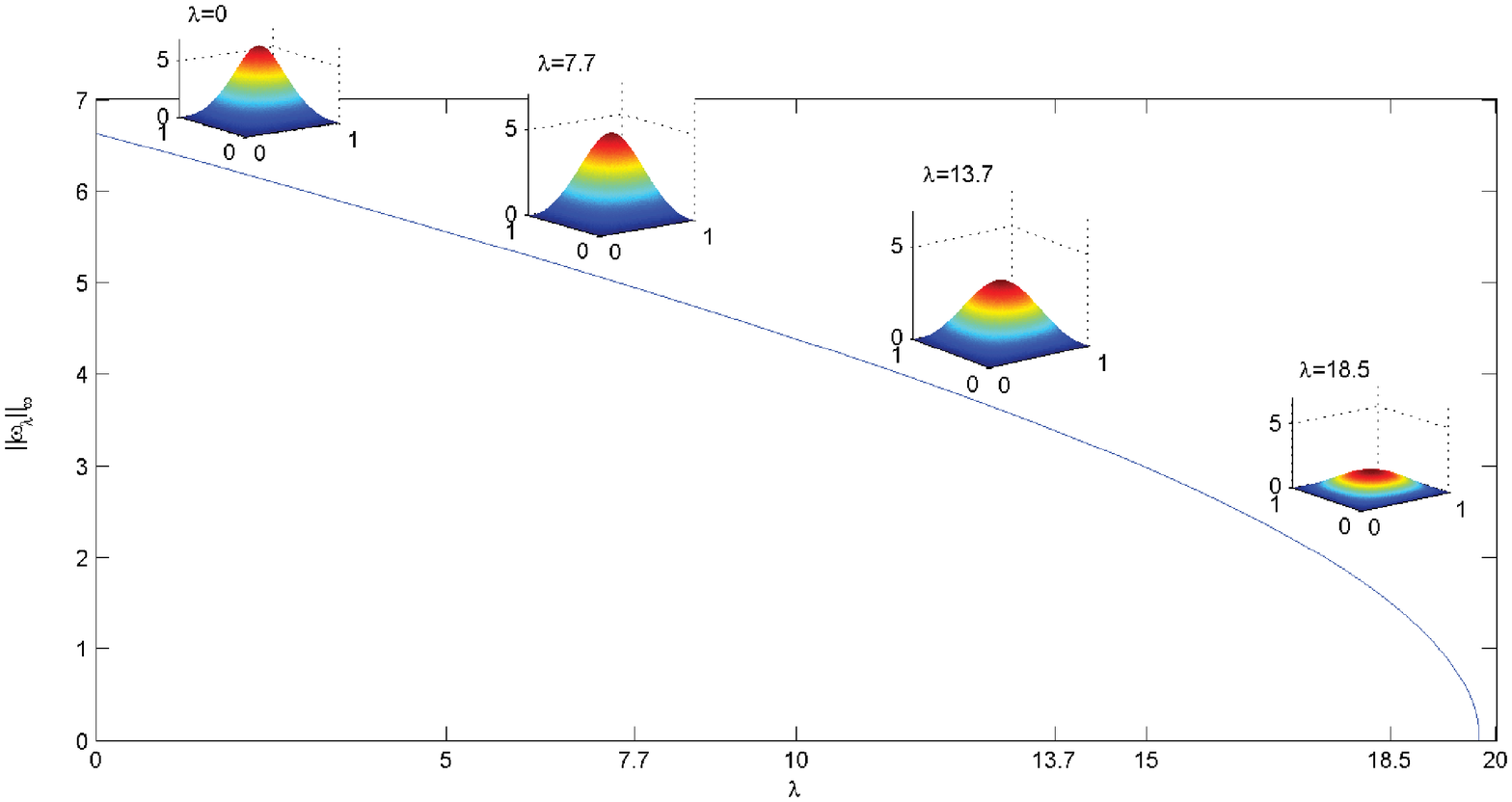}
	\label{fig:wmaxplotp3}
	\caption{Curve $(\lambda,\|\omega_\lambda\|_\infty)$ with samples of $\omega_\lambda$ in the
case $p=3$}
\end{figure}

Using $\bar{\lambda} = 18.5$ (which is not the minimally possible choice; e.g. $\bar\lambda =15.7$ could have been chosen) and $M+1=94$ values
$0=\lambda^0 < \lambda^1 < \cdots < \lambda^{93} = 18.5$ (with $\lambda^{1}=0.1, \lambda^{2}=0.3$
and the remaining gridpoints equally spaced with distance 0.2) we computed approximations $\omega^0,
\dots, \omega^{93}$ with $N=16$ in \eqref{eq6.2}, as well as defect bounds $\delta^0,
\dots,\delta^{93}$ and constants $K^0, \dots, K^{93}$, by the methods described in Section 6 and 7.

Figure 1 shows an approximate branch $[0,2 \pi^2) \to \R, \lambda \mapsto
\|\omega_{\lambda}\|_\infty$. The continuous plot has been created by interpolation of the above grid
points $\lambda^j$, plus some more grid points between $18.5$ and $2\pi^2$, where we computed additional approximations.

For some selected values of $\lambda$, Table 1 shows, with an obvious sub- and superscript notation
for enclosing intervals, the computed eigenvalue bounds for problem \eqref{eq7.6} (giving $K_{\lambda}$ by
\eqref{eq7.7}). These were obtained using the Rayleigh-Ritz and the Lehmann-Goerisch method, and the
homotopy method briefly mentioned at the end of Section 7 (exploiting also the symmetry considerations addressed in Remark 7.1).
The integer $m$, needed for these procedures, has been chosen different
(between $3$ and $10$) for different values of $\lambda$, according to the outcome of the homotopy. This resulted in a slightly different
quality of the eigenvalue enclosures. 
{
\renewcommand{\arraystretch}{1.5}
\begin{table}[H]
\centering
\begin{tabular}{|c|l|l|}
\hline
&$\kappa_1$ & $\kappa_2$\\
\hline\hline
$\omega_0 $ &  $0.34350814513_{229}^{840}$& $2.492570_{450}^{712}$	\\\hline
$\omega_{2.7}$	&	$0.37521912233_{290}^{850}$	& $2.6221837_{393}^{653}$	\\\hline
$\omega_{6.7}$ &  $0.4373273950_{355}^{411}$ & $2.87378161_{204}^{409}$  \\\hline
$\omega_{10.7}$ &  $0.52752354636_{169}^{621}$	& $3.223417042_{185}^{515}$	\\\hline
$\omega_{14.7}$ &	$0.6676848259_{379}^{417}$	& $3.725209290_{830}^{988}$ \\\hline
$\omega_{18.5}$ & 	$0.89237445994_{555}^{742}$ & $4.46288110_{093}^{102}$ \\
\hline
\end{tabular}
\vspace{0.2cm}
\caption{Eigenvalue enclosures for the first two eigenvalues}
\end{table}
}

%

Table 2 contains, for some selected of the 186 $\lambda$-half-intervals,
\begin{enumerate}
\item[a)] the defect bounds $\delta_{\lambda}$ obtained by \eqref{eq8.7} from the grid-point defect
bounds $\delta^{i-1}, \delta^i$, and from the grid-width characteristics $\rho_i, \tau_i$ defined in
\eqref{eq8.5}, \eqref{eq8.6},
\item[b)] the constants $K_{\lambda}$ obtained by \eqref{eq8.10} from the grid-point constants $K^i$
and the grid-width parameters $\nu_i$ defined in \eqref{eq8.8} (note that $\mu_i = 0.1$ for all
$i$),
\item[c)] the error bounds $\alpha_{\lambda}$ computed according to \eqref{eq3.4}, \eqref{eq3.5}.
\end{enumerate}

Thus, Corolllary \ref{coro2.1}, together with all the considerations in the previous sections,
proves Theorem \ref{theo1.1}.
\begin{table}[H]
\centering
\begin{tabular}{|c|c|c|c|}
\hline
$\lambda$-interval	&   $\delta_\lambda$	&	$K_\lambda$		&	$\alpha_\lambda$ \\
\hline \hline
[0,0.05)			&		0.0005943		&	1.7443526	&	0.0010378		\\
\hline							
(2,2.1)			&	0.0023344		&	1.7707941		&	0.0041521	\\
\hline
(6,6.1)			&	0.0022937	&	1.6669879	&	 0.0038369		\\
\hline
(10,10.1) 	&	0.0023644	&	1.5677657	&	0.0037168	\\
\hline
(14,14.1) 	&	0.0026980	&	1.9582604	&	0.0053028\\
\hline
(16,16.1)	&	0.0031531	&	3.2267762	& 	0.0102701	\\
\hline
(18.4,18.5] 	&	0.0050056  &	13.8930543	&	0.0882899	\\
\hline
\end{tabular}
\vspace{0.2cm}
\caption{}
\end{table}

{\bf Acknowledgments.} The authors are grateful to two anonymous referees for their helpful remarks and suggestions.
 


\vspace{1cm}

\bibliographystyle{siam}
\bibliography{bibliography}

\begin{thebibliography}{10}

\bibitem{AdiYad94}
{\sc Adimurthi and S.~Yadava}, {\em An elementary proof for the uniqueness of
  positive radial solution of a quasilinear {D}irichlet problem}, Arch. Rat.
  Mech. Anal., 126 (1994), pp.~219--229.

\bibitem{AftPac03}
{\sc A.~Aftalion and F.~Pacella}, {\em Uniqueness and nondegeneracy for some
  nonlinear elliptic problems in a ball}, Journ. Diff. Eq., 195 (2003),
  pp.~380--397.

\bibitem{Beh87}
{\sc H.~Behnke}, {\em Inclusion of eigenvalues of general eigenvalue problems
  for matrices}, in: U. Kulisch, H. J. Stetter (Eds.), Scientific Computation
  with Automatic Result Verification, Computing 6 (Suppl.) (1987), pp.~69--78.

\bibitem{BehGoe94}
{\sc H.~Behnke and F.~Goerisch}, {\em Inclusions for eigenvalues of selfadjoint
  problems}, in: J. Herzberger (Ed.), Topics in Validated Computations, Series
  Studies in Computational Mathematics, North-Holland, Amsterdam (1994),
  pp.~277--322.

\bibitem{BreuHorMcKennaPlum06}
{\sc B.~Breuer, J.~Horak, P.~J. McKenna, and M.~Plum}, {\em A computer-assisted
  existence and multiplicity proof for travelling waves in a nonlinearly
  supported beam}, J. Differential Equations, 224 (2006), pp.~60--97.

\bibitem{BreuMcKennaPlum03}
{\sc B.~Breuer, P.~J. McKenna, and M.~Plum}, {\em Multiple solutions for a
  semilinear boundary value problem: a computational multiplicity proof}, J.
  Differential Equations, 195 (2003), pp.~243--269.

\bibitem{CR1}
{\sc M.~G. Crandall and P.~H. Rabinowitz}, {\em Bifurcation from simple
  eigenvalues}, J. Functional Analysis, 8 (1971), pp.~321--340.

\bibitem{CR2}
\leavevmode\vrule height 2pt depth -1.6pt width 23pt, {\em Bifurcation,
  perturbation of simple eigenvalues and linearized stability}, Arch. Rational
  Mech. Anal., 52 (1973), pp.~161--180.

\bibitem{DamGroPac99}
{\sc L.~Damascelli, M.~Grossi, and F.~Pacella}, {\em Qualitative properties of
  positive solutions of semilinear elliptic equations in symmetric domains via
  the maximum principle}, Ann. Inst. H. Poincar\'e, 16 (1999), pp.~631--652.

\bibitem{Dan88}
{\sc E.~N. Dancer}, {\em The effect of the domain shape on the number of
  positive solutions of certain nonlinear equations}, Journ. Diff. Eq., 74
  (1988), pp.~120--156.

\bibitem{GidNiNir79}
{\sc B.~Gidas, W.~M. Ni, and L.~Nirenberg}, {\em Symmetry and related
  properties via the maximum principle}, Comm. Math. Phys., 68 (1979),
  pp.~209--243.

\bibitem{Gro00}
{\sc M.~Grossi}, {\em A uniqueness result for a semilinear elliptic equation in
  symmetric domains}, Adv. Diff. Equations, 5 (2000), pp.~193--212.

\bibitem{Kla93}
{\sc R.~Klatte, U.~Kulisch, C.~Lawo, M.~Rausch, and A.~Wiethoff}, {\em C-XSC-A
  C++ Class Library for Extended Scientific Computing}, Springer, Berlin, 1993.

\bibitem{Lad68}
{\sc O.~A. Ladyzhenskaya and N.~N. Ural'tseva}, {\em Linear and quasilinear
  elliptic equations}, Academic Press, New York-London, 1968.

\bibitem{Leh63}
{\sc N.~J. Lehmann}, {\em Optimale {E}igenwerteinschlie{\ss}ungen}, Numer.
  Math., 5 (1963), pp.~246--272.

\bibitem{McKPacPlR09}
{\sc P.~J. McKenna, F.~Pacella, M.~Plum, and D.~Roth}, {\em A uniqueness result
  for a semilinear elliptic problem: A computer-assisted proof}, J.
  Differential Equations, 247 (2009), pp.~2140--2162.

\bibitem{NagNaYa99}
{\sc K.~Nagatou, M.~T. Nakao, and N.~Yamamoto}, {\em An approach to the
  numerical verification of solutions for nonlinear elliptic problems with
  local uniqueness}, Numer. Funct. Anal. Optim., 20 (1999), pp.~543--565.

\bibitem{NaYa95}
{\sc M.~T. Nakao and N.~Yamamoto}, {\em Numerical verifications for solutions
  to elliptic equations using residual iterations with higher order finite
  elements}, J. Comput. Appl. Math., 60 (1995), pp.~271--279.

\bibitem{NiNuss85}
{\sc W.~M. Ni and R.~D. Nussbaum}, {\em Uniqueness and nonuniqueness for
  positive radial solutions of ${\Delta} u + f(u,\tau) = 0$}, Comm. Pure Appl.
  Math., 38 (1985), pp.~67--108.

\bibitem{PacSri03}
{\sc F.~Pacella and P.~N. Srikanth}, {\em Solutions of semilinear problems in
  symmetric planar domains, {ODE} behaviour and uniqueness of branches},
  Progress in Nonlinear Diff. Eq. and Their Appl., 54 (2003), pp.~239--244.

\bibitem{Plum92}
{\sc M.~Plum}, {\em Explicit ${H}^2$-estimates and pointwise bounds for
  solutions of second-order elliptic boundary value problems}, J. Math. Anal.
  Appl., 165 (1992), pp.~36--61.

\bibitem{Plum97}
\leavevmode\vrule height 2pt depth -1.6pt width 23pt, {\em Guaranteed numerical
  bounds for eigenvalues}, in: D. Hinton, P. W. Schaefer (Eds.), Spectral
  Theory and Computational Methods of Sturm-Liouville Problems, Marcel Dekker,
  New York (1997), pp.~313--332.

\bibitem{Plum08}
\leavevmode\vrule height 2pt depth -1.6pt width 23pt, {\em Existence and
  multiplicity proofs for semilinear elliptic boundary value problems by
  computer assistance}, DMV Jahresbericht, JB 110 (2008), pp.~19--54.

\bibitem{PlumWie02}
{\sc M.~Plum and C.~Wieners}, {\em New solutions of the {G}elfand problem}, J.
  Math. Anal. Appl.,  (2002), pp.~588--606.

\bibitem{Rab71}
{\sc P.~H. Rabinowitz}, {\em Some global results for nonlinear eigenvalue
  problems}, Journ. Funct. Anal., 7 (1971), pp.~487--513.

\bibitem{Rek80}
{\sc K.~Rektorys}, {\em Variational Methods in Mathematics, Science and
  Engineering}, Reidel Publ. Co., Dordrecht, 1980.

\bibitem{Rump02}
{\sc S.~M. Rump}, {\em INTLAB-INTerval LABoratory, a Matlab toolbox for
  verified computations, Version 4.2.1, Inst. Informatik, TU Hamburg-Harburg,
  2002}, http://www.ti3.tu-harburg.de/rump/intlab/.

\bibitem{Sri93}
{\sc P.~N. Srikanth}, {\em Uniqueness of solutions of nonlinear {D}irichlet
  problems}, Diff. Int. Eq., 6 (1993), pp.~663--670.

\bibitem{Zhang92}
{\sc L.~Zhang}, {\em Uniqueness of positive solutions of ${\Delta} u + u^p + u
  = 0$ in a finite ball}, Comm. Part. Diff. Eq., 17 (1992), pp.~1141--1164.

\bibitem{Zou94}
{\sc H.~Zou}, {\em On the effect of the domain geometry on the uniqueness of
  positive solutions of ${\Delta} u + u^p = 0$}, Ann. Sc. Norm. Sup. Pisa, 3
  (1994), pp.~343--356.

\end{thebibliography}

\end{document}